\documentclass[12pt]{amsart}
\usepackage{amsfonts, amsmath, latexsym, epsfig}
\usepackage{amssymb}
\usepackage{epsf}
\usepackage{url}

\newcommand{\RR}{\ensuremath{\mathbb{R}}}

\newcommand{\QQ}{\ensuremath{\mathbb{Q}}}

\newcommand{\ZZ}{\ensuremath{\mathbb{Z}}}

\newtheorem{theorem}{Theorem}

\def\QuotS#1#2{\leavevmode\kern-.0em\raise.2ex\hbox{$#1$}\kern-.1em/\kern-.1em\lower.25ex\hbox{$#2$}}

\DeclareMathOperator{\conv}{conv}

\begin{document}

\author{Mathieu Dutour Sikiri\'c}
\address{Mathieu Dutour Sikiri\'c, Rudjer Boskovi\'c Institute, Bijenicka 54, 10000 Zagreb, Croatia, Fax: +385-1-468-0245}
\email{mathieu.dutour@gmail.com}

\title{The Birkhoff polytope of the groups $\mathsf{F}_4$ and $\mathsf{H}_4$}
\date{}

\maketitle

\begin{abstract}
We compute the set of facets of the polytope which is the convex hull
of the Coxeter groups $\mathsf{F}_4$ or $\mathsf{H}_4$:
\begin{itemize}
\item For the group $\mathsf{F}_4$ we found $2$ orbits of facets which
  contradicts previous results published in \cite{birkhoff}.
\item For the group $\mathsf{H}_4$ we found $1063$ orbits of facets which
  provides a counterexample to the conjecture of \cite{birkhoff}.
\end{itemize}
\end{abstract}

\section{Introduction}
Given a finite group $G$ acting linearly on a real vector space $\RR^n$, there
is a strong interest in finding the orbits of facets of the convex
hull $Gx$ of a vector $x$.

For a coxeter group $G$ with its natural action on $\RR^n$, the structure
is well known and given by the Wythoff construction (See \cite{Wythoff}).
Other representations were considered in \cite{Alternahedron} for the
alternating group. Another very interesting case for a group having a
$n$-dimensional representation is to consider the action of the group
on itself. For the symmetric group $S_n$ this gets us the Birkhoff
polytope.

In \cite{ConvexHullsCoxeterGroup,birkhoff} the description was extended
to other Coxeter groups by introducing the Birkhoff tensors $B_G$
which are facets of $\conv\,G$. In \cite{zobin} the symmetry group of
$\conv\, G$ are determined for all the finite Coxeter groups.

The authors proved the following result:
\begin{theorem}
For $G=\mathsf{A}_n$, $\mathsf{B}_n$, $\mathsf{I}_2(n)$ and $\mathsf{H}_3$ all facets of
$\conv\, G$ are Birkhoff tensors.
\end{theorem}
They also proved that for $\mathsf{D}_4$ the result does not hold and
they claim in Theorem 8.1 that this implies that the result does not
hold for $\mathsf{D}_n$ and $\mathsf{E}_n$.
Note that the authors also claimed that $\mathsf{F}_4$ satisfy the theorem but
we prove that this is not true.

The authors conjectured \cite[Problem 8.1]{birkhoff} that for $\mathsf{H}_4$ the
result does hold.
As it turns out, this is not true since while there is just
one orbit of Birkhoff tensors, there is more than one orbit of facets:

\begin{theorem}
The polytope $\conv\, G$ has $188455824000$ facets in $1063$ orbits.
\end{theorem}

The proof of this result is computational with the algorithms presented in
Section \ref{SEC_algorithm} and the results presented in Section \ref{SEC_results}.

\section{Algorithms}\label{SEC_algorithm}

The effective computation of dual description of polytopes is a classic
problem in \cite{HandbookDiscCompGeometry}. We designed software for
using symmetries of polytope when computing their dual descriptions.
The initial version of the code was written in {\tt GAP} with some parts
in {\tt C}/{\tt C++}. After several extensions, the code was completely
ported into {\tt C++}. As a side product of that, we are allowed to
select the numeric type of the occurring matrix entries.

\subsection{Dual description algorithms}

In order to compute the facets of the polytope, we used the recursive adjacency
decomposition technique. This method has been used for many different computations
and is explained in \cite{perfectdim8,montreal,ContactLeech,CUTsmallGraphs}.

The idea of the adjacency decomposition technique is to compute one facet and from
this facet to compute the adjacent facets. The obtained facets are then tested
for equivalence. Computing the facets adjacent to a given facet is itself a
dual description problem, therefore one may need to apply the method recursively
hence the name recursive adjacency decomposition method.

Some early termination criterion are given in \cite{CUTsmallGraphs} and allow us
to avoid having to compute the adjacencies of all the facets.

For a long time we used the code developed in \cite{polyhedral_gap} which is a
package of \cite{gap}. For several reasons we have developed a new {\tt C++}
implementation (see \cite{polyhedral_cpp}) that allow us to gain additional speed
and functionality.

\subsection{Fields}

The commonly used numerical type is {\bf mpq\_class} from the GMP library
which a multiprecision rational type. It is supported in GAP as well as {\tt C++}.
In order to compute with polytopes related to the Coxeter groups $\mathsf{H}_3$ and $\mathsf{H}_4$
one needs to allow for the ring $\QQ\lbrack \sqrt{5}\rbrack$.
Implementing the arithmetic operations ($+$, $-$, $*$, $/$) is relatively easy
but the sign determinations require more care.

Testing if $a + b \sqrt{5}$ is positive can be done in the following way.
If $a$ and $b$ are of the same sign it is easy to conclude.
Otherwise, $a$ and $b$ are of opposite sign and we write
\begin{equation*}
a + b \sqrt{5} = \frac{a - b \sqrt{5}}{a^2 - 5 b^2}
\end{equation*}
The sign of $a - b \sqrt{5}$ can be decided and together with the sign of the
rational number $a^2 - 5 b^2$ we can conclude.
The same strategy allows one to decide the sign in mixed cubic rings.
For other real fields of algebraic numbers, different approaches would have
to be used.

Also, for some subroutines like {\tt lrs} one needs only to use ring operations.
In that case one can reduce to the case of $\ZZ\lbrack \sqrt{5}\rbrack$. For the kind
of computations we are doing here there is no need for algebraic closures or such
kind of constructions.

\subsection{Canonicalization strategies}

In preceding works, when we had two orbits of facets in order to check isomorphism,
we used the \cite{gap} implementation of the partition backtrack. See \cite{Leon1,Leon2}
and \cite{JeffersonReExplaination} for accounts of this class of algorithms.
That is we encode facets by the subset of their incident vertices and then use the
partition backtrack for set equivalence.

However, in \cite{MinimalCanonicalImage} an algorithm for finding a canonical
representative of a subset for a permutation group action was found. This greatly
simplifies the code since for $N$ orbits instead of having to compute $N$ equivalences,
we simply have to do one canonicalization and one string comparison.

\section{Results}\label{SEC_results}

\subsection{Results for $\mathsf{F}_4$}
For the group $\mathsf{F}_4$ we use the following two generators:
\begin{equation*}
\frac{1}{2}\left(\begin{array}{cccc}
1 & -1 & 1 & 1\\
-1 & -1 & -1 & 1\\
1 & 1 & -1 & 1\\
1 & -1 & -1 & -1
\end{array}\right)
,
\frac{1}{2}\left(\begin{array}{cccc}
-1 & 1 & -1 & -1\\
1 & -1 & -1 & -1\\
-1 & -1 & 1 & -1\\
1 & 1 & 1 & -1
\end{array}\right)
\end{equation*}

There are exactly two orbits of facets for $\conv(\mathsf{F}_4)$:
\begin{enumerate}
\item Orbit 1 of Birkhoff tensor with incidence $288$. Stabilizer has order $4608$. One representative inequality is $Tr(XA) \leq 1$ with
\begin{equation*}
A = \left(\begin{array}{cccc}
0 & 0 & 0 & 0\\
0 & 0 & 0 & 0\\
0 & 0 & 0 & 0\\
1 & 0 & 0 & -1
\end{array}\right)
\end{equation*}
\item Orbit 2 of facets with incidence $288$. Stabilizer has order $48$. One representative inequality is $Tr(XA) \leq 1$ with
\begin{equation*}
A = \frac{1}{4}\left(\begin{array}{cccc}
1 & 0 & 1 & 0\\
0 & 1 & 0 & -1\\
0 & 1 & 0 & 1\\
1 & 0 & 1 & -2
\end{array}\right)
\end{equation*}
\end{enumerate}
The total number of facets is $55872$.

\subsection{Results for $\mathsf{H}_4$}
For the group $\mathsf{H}_4$ we use the following two generators:
\begin{equation*}
\frac{1}{4}\left(\begin{array}{cccc}
1 & -2 & -1 & 0\\
2 & 2 & -2 & 2\\
1 & -2 & 0 & 1\\
0 & -2 & -1 & 1
\end{array}\right)
+\frac{\sqrt{5} }{4}\left(\begin{array}{cccc}
1 & 0 & 1 & 0\\
0 & 0 & 0 & 0\\
-1 & 0 & 0 & 1\\
0 & 0 & -1 & -1
\end{array}\right)
,\left(\begin{array}{cccc}
-1 & 0 & 0 & 0\\
0 & 0 & 0 & -1\\
0 & -1 & 0 & 0\\
0 & 0 & 1 & 0
\end{array}\right)
\end{equation*}

We found $1063$ orbits of facets of $\conv\,\mathsf{H}_4$.
If we express such facets in the form $Tr(AX)\leq 1$ then one $1$, $4$, $130$, $928$ orbits of facets of rank $1$, $2$, $3$ and $4$.

\begin{table}
\begin{center}
\begin{tabular}{|cc|cc|cc|cc|}
\hline
$p$  & Nr.      &    $p$  & Nr.      &    $p$  & Nr.      &    $p$  & Nr.\\\hline
$16$ & $376$    &    $17$ & $282$    &    $18$ & $116$    &    $19$ & $85$\\
$20$ & $48$    &    $21$ & $30$    &    $22$ & $12$    &    $23$ & $5$\\
$24$ & $44$    &    $25$ & $3$    &    $26$ & $31$    &    $28$ & $10$\\
$30$ & $7$    &    $32$ & $5$    &    $36$ & $4$    &    $38$ & $1$\\
$48$ & $1$    &    $100$ & $1$    &    $120$ & $1$    &    $480$ & $1$\\
\hline
\end{tabular}
\caption{For each incidence $p$ the number of orbits of incidence $p$ is given.}\label{Tab_Incidence}
\end{center}
\end{table}

\begin{table}
\begin{center}
\begin{tabular}{|cc|cc|cc|cc|}
\hline
$s$  & Nr.      &    $s$  & Nr.      &    $s$  & Nr.      &    $s$  & Nr.\\\hline
$1$ & $800$    &    $2$ & $189$    &    $4$ & $50$    &    $6$ & $1$\\
$8$ & $8$    &    $12$ & $4$    &    $16$ & $3$    &    $24$ & $1$\\
$36$ & $1$    &    $40$ & $1$    &    $48$ & $2$    &    $120$ & $1$\\
$576$ & $1$    &    $2880$ & $1$    &       &          &         &          \\
\hline
\end{tabular}
\caption{For each size $s$ the number of orbits of orbits having a stabilizer of size $s$ is given.}\label{Tab_StabSiz}
\end{center}
\end{table}

Tables \ref{Tab_Incidence} and \ref{Tab_StabSiz} give the statistics about the incidence of orbits and about the size of their
stabilizers. The full list of orbits is presented in \cite{convcoxeterresult}. See below one matrix of incidence $120$ with a
stabilizer of size $120$. This suffices to show that the conjecture is false.
\begin{equation*}
A=\frac{1}{4}\left(\begin{array}{cccc}
-6 & -11 & -7 & 4\\
0 & -2 & 3 & 1\\
0 & -1 & 4 & 3\\
0 & 0 & 0 & 0
\end{array}\right)
+\frac{\sqrt{5} }{4}\left(\begin{array}{cccc}
2 & 5 & 3 & -2\\
0 & 0 & -1 & -1\\
0 & 1 & -2 & -1\\
0 & 0 & 0 & 0
\end{array}\right)
\end{equation*}

\providecommand{\bysame}{\leavevmode\hbox to3em{\hrulefill}\thinspace}
\providecommand{\MR}{\relax\ifhmode\unskip\space\fi MR }
\providecommand{\MRhref}[2]{%
  \href{http://www.ams.org/mathscinet-getitem?mr=#1}{#2}
}
\providecommand{\href}[2]{#2}

\end{document}